\documentclass[10pt,twocolumn,letterpaper]{article}

\usepackage{cvpr}
\usepackage{graphicx}
\usepackage{comment}
\usepackage{amsmath,amssymb} 
\usepackage{color}
\usepackage{times}
\usepackage{epsfig}
\usepackage{diagbox}
\usepackage{afterpage}
\usepackage{makecell}
\usepackage{lineno}
\usepackage{amsfonts}       
\usepackage{amsthm}   

\usepackage{bm}
\usepackage{algorithm}
\usepackage{algorithmic}

\newcommand{\be}{\begin{eqnarray}}
\newcommand{\ee}{\end{eqnarray}}
\newcommand{\beq}{\begin{equation}\begin{aligned}}
\newcommand{\eeq}{\end{aligned}\end{equation}}
\newcommand{\beqn}{\begin{equation*}\begin{aligned}}
\newcommand{\eeqn}{\end{aligned}\end{equation*}}
\newcommand{\ben}{\begin{eqnarray*}}
\newcommand{\een}{\end{eqnarray*}}

\usepackage{algorithm}
\usepackage{algorithmic}

\usepackage{multirow}
\usepackage{booktabs}       

\usepackage[pagebackref=true,breaklinks=true,letterpaper=true,colorlinks,bookmarks=false]{hyperref}

\usepackage{authblk}
\cvprfinalcopy
\ifcvprfinal\fi
\begin{document}

\title{Deep Domain Decomposition Method: Elliptic Problems} 
\author[1,2]{Wuyang Li}
\author[1]{Xueshuang Xiang\thanks{Corresponding author: xiangxueshuang@qxslab.cn}}
\author[2]{Yingxiang Xu}
\affil[1]{\normalsize Qian Xuesen Laboratory of Space Technology, China Academy of Space Technology}
\affil[2]{\normalsize School of Mathematics and Statistics, Northeast Normal University}
\date{}

\maketitle

\begin{abstract}
This paper proposes a deep-learning-based domain decomposition method (DeepDDM), which leverages deep neural networks (DNN) to discretize the subproblems divided by domain decomposition methods (DDM) for solving partial differential equations (PDE).
Using DNN to solve PDE is a physics-informed learning problem with the objective involving two terms, domain term and boundary term, which respectively make the desired solution satisfy the PDE and corresponding boundary conditions. 
DeepDDM will exchange the subproblem information across the interface in DDM by adjusting the boundary term for solving each subproblem by DNN. 
Benefiting from the simple implementation and mesh-free strategy of using DNN for PDE, DeepDDM will simplify the implementation of DDM and make DDM more flexible for complex PDE, e.g., those with complex interfaces in the computational domain. 
This paper will firstly investigate the performance of using DeepDDM for elliptic problems, including a model problem and an interface problem.
The numerical examples demonstrate that DeepDDM exhibits behaviors consistent with conventional DDM: the number of iterations by DeepDDM is independent of network architecture and decreases with increasing overlapping size. 
The performance of DeepDDM on elliptic problems will encourage us to further investigate its performance for other kinds of PDE and may provide new insights for improving the PDE solver by deep learning.
\end{abstract}

\section{Introduction}
In the realms of engineering and scientific computing, the domain decomposition methods (DDM) has attracted a great deal of interest and is well-known as an efficient approach for solving partial differential equations (PDE).
The main idea of DDM is first splitting the computational domain into smaller subdomains and then solving in parallel the subproblems defined in subdomains, along with exchanging solution information between adjacent subdomains. 
When combined with the discretization of PDE by finite element methods (FEM) or finite difference methods (FDM), DDM can achieve remarkable performance. 
This paper will propose an approach named DeepDDM which investigates the behavior of discretizing PDE by deep neural networks (DNN) in DDM. 
Compared with DDM-FEM/DDM-FDM, DeepDDM will simplify the implementation of DDM. 
Especially, DeepDDM is more flexible for complex PDE, e.g., those with complex interfaces in the computational domain.  

Due to the success of deep learning (DL) in engineering, the research using deep learning to solve PDE has provoked interest from numerous researchers. 
Recently, many relevant papers have been published focusing on different topics.
\cite{raissi2019physics} presents physics-informed neural networks (PINN) that take initial conditions and boundary conditions as the penalties of the optimization objective loss function. 
By using the backpropagation algorithm and automatic differentiation, it realizes the calculation of the high-order differential in the framework of TensorFlow.
The network can achieve good accuracy for both forward and inverse problems.
\cite{sirignano2018dgm} proposes a deep Galerkin method, which is similar to PINN in that initial conditions and boundary conditions are also added as penalty terms into the optimization objective loss function.
The difference is that \cite{sirignano2018dgm} bypasses the high-order differentiation of the neural network through the Monte Carlo method.
For well-known reasons, the computational cost of high-order differentiation of neural networks is exorbitant.
Another method, a deep Ritz method, is presented by \cite{weinan2017deep}, in which different network architectures and loss functions are designed for variational problems, particularly the ones that arise from PDE.
The large multiple layer convolutional neural network is used to solve and discover evolutionary PDE by \cite{long2019pde}.
Recently, a large number of researchers have attempted to use deep learning tools to solve various practical problems \cite{berg2018unified,rudd2013constrained}.

In the field of combining machine learning and domain decomposition, \cite{mai2002mesh} presented a combination of radial basis function network methods and domain decomposition technique for approximating functions and solving Poisson equations. 
Recently, as we prepared for this paper, some related works of combining domain decomposition method and deep learning have appeared on the Internet. 
\cite{li2019d3m} presented D3M, a method that combines Deep Ritz method and DDM to solve general PDE in parallel. 
However, our approach focus on the performance of combing PINN and DDM for Elliptic problems, especially for the problem with a complex interface. 
\cite{dwivedi2019distributed} presented a Distributed PINN that divides the problem into many subproblems according to a domain splitting, and then use the total loss, i.e. the sum of all subproblem’s loss and the loss on the interface, to train all the subproblems at the same time. We notice that the training objective of each subproblem in Distributed PINN is related to its neighbors at each training step. 
That means Distributed PINN is not DDM-style method, since we need to first independently solve each subproblem and then exchange the interface information in DDM, like D3M and the proposed DeepDDM. 
On the application front, \cite{kissas2020machine} applied PINN in different domains of cardiovascular flows modeling with the similar spirit of Distributed PINN. 

We present a deep-learning-based domain decomposition framework, called the deep domain decomposition method (DeepDDM), which extracts the spirit of deep learning and domain decomposition. 
Using DNN to solve PDE is a physics-informed learning problem with the objective involving two terms, domain term and boundary term, which respectively make the desired solution satisfy the PDE and corresponding boundary conditions \cite{raissi2019physics}.
By dividing the domain of interest into many subdomains, DeepDDM alternatively solves each subproblem by DNN and exchanges interface information until convergence. 
The decomposition of the domain of interest is based on the properties of the original physical background or the convenience of computing. 
Since DNN makes solving PDE a simple learning problem and is actually a mesh-free strategy, DeepDDM will simplify the implementation of DDM and make DDM more flexible for complex PDE. 

This paper will first investigate the performance of DeepDDM on elliptic problems, including a model problem and an interface problem. 
For the model problem with the Dirichlet boundary condition, we consider the 
original problem divided into two or four subproblems in coordination with overlapping information and varying discrete functional spaces (network architectures).
The results demonstrate that DeepDDM has properties similar to DDM with FEM or FDM:
the necessary number of iterations is independent of the network architecture, given the number of subproblems; 
the necessary number of iterations increases along with the number of subproblems; 
the necessary number of iterations decreases according to increasing overlap size; and
the numerical convergence rate coincides with the analytic convergence rate in various cases. 
Next, we investigate the performance of DeepDDM on an artificial interface problem, which is divided into two subproblems. 
Although only the elementary Dirichlet-Neumann interface condition is adopted, DeepDDM also exhibits good performance for the interface problem, just as for the model problem.
In addition, DeepDDM reaches the relative error accuracy of $10^{-3}$ within several iterations for varying network architectures, especially for the case that the coefficient contrast (the ratio of the maximal diffusion coefficient on the minimal diffusion coefficient of the interface problem) is $20$. 
In conclusion, as an approach with simple implementation, DeepDDM produces similar results as DDM with FEM or FDM and may improve the performance of using DDM for interface problems. 

This paper is structured as follows. 
In \autoref{sec:Schwarz} and \autoref{sec:DeepLearning}, we provide a brief review of DDM and
the PINN, respectively. 
A detailed description of DeepDDM for the general differential operator and boundary condition is presented in \autoref{sec:DeepDDM}. 
We implement and test the performance of DeepDDM on two kinds of elliptic problems in \autoref{sec:Examples}. 
The \autoref{sec:Conclusion} will conclude this paper and present some directions for future work. 

\begin{figure}[!t]
	\centering
	\includegraphics[width=0.4\textwidth, height=\textheight, keepaspectratio]
	{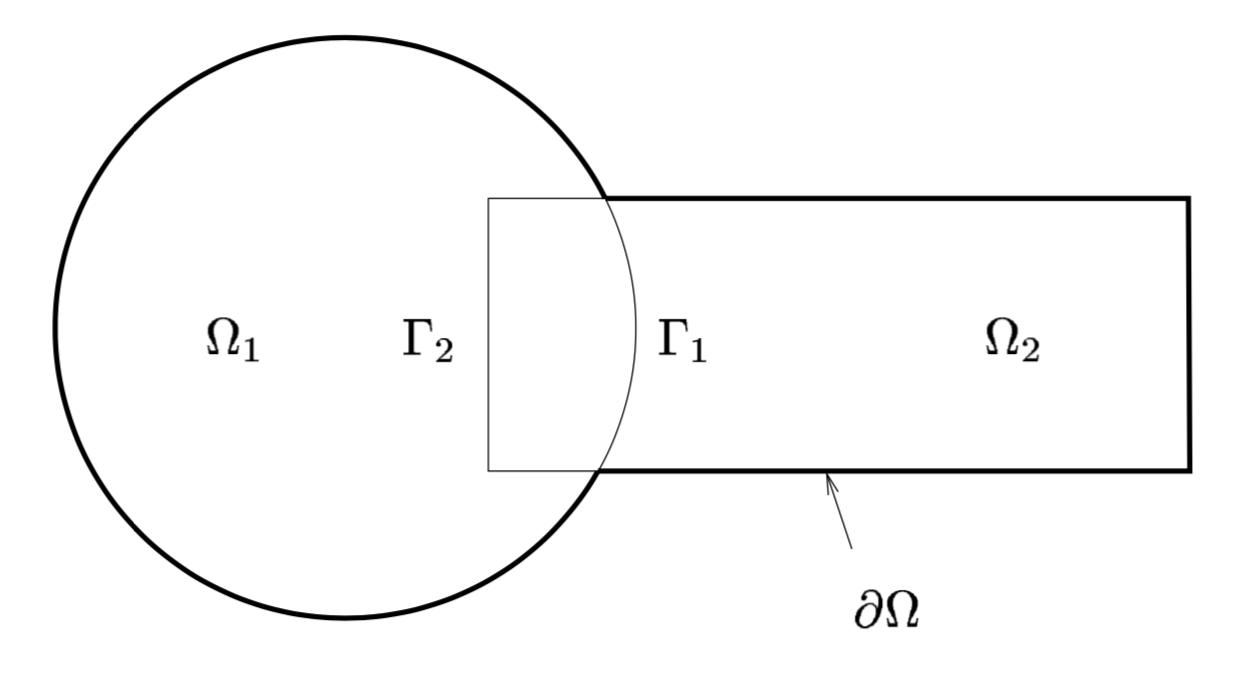}
	\caption{The diagram used by Schwarz in 1870.}
	\label{fig:diagram_schwarz}
\end{figure}

\section{Domain decomposition methods}\label{sec:Schwarz}
DDM are parallel, potentially fast, robust algorithms for solving PDE, which has been discretized using, e.g., FEM or FDM. 
Even in the 19th century, Hermann Schwarz considered a Poisson problem set on a union of simple geometries and introduced an alternating iterative method in \cite{schwarz1870ueber}.
Nearly a century later, Pierre-Louis Lions presented parallel Schwarz methods for parallel computing in \cite{lions1988schwarz} after the parallel computing capability became available.
For the discrete system, the multiplicative Schwarz method, sequential algorithms in nature, is introduced in \cite{chan1994domain,smith2004domain} and the additive Schwarz method, parallel algorithms, is studied by Dryja and Widlund in \cite{dryja1990some}.
With the wide use of parallel computers, DDM develop rapidly and provide many fast iterative algorithms and efficient preconditioners.
Balancing domain decomposition by constraints (BDDC) and finite element tearing and interconnecting (FETI) have been researched for a variety of PDE as typical representatives of primal methods and dual methods, respectively~\cite{dolean2015introduction,lanser2015nonlinear,mathew2008domain,pechstein2017unified}.
Moreover, a hybrid method between primal and dual, dual-primal finite element tearing and interconnecting (FETI-DP), enforcing equality of the solution at the subdomain interface by Lagrange multipliers, except at the subdomain corner, was initially introduced by \cite{farhat2001feti} and has developed into many varieties~\cite{dolean2015introduction,klawonn2014nonlinear}.

Here we briefly introduce a classical DDM, the Schwarz alternating 	method of Schwarz in 1870, as follows,
\beq
\Delta u &= 0, \quad \text{in}\,\Omega, \\
u &= g, \quad \text{on}\,\partial \Omega,
\eeq
where $\Omega$, shown as \autoref{fig:diagram_schwarz}, may be divided into a disk $\Omega_1$ and a rectangle $\Omega_2 $, $\Omega=\Omega_1 \cup \Omega_2$, with interfaces $\Gamma_1:=\partial\Omega_1 \backslash \partial\Omega $ and $\Gamma_2:= \partial\Omega_2 \backslash \partial\Omega $, and $ g $ is a function given as the boundary condition. 
The Schwarz alternating method starts with an initial guess $u_2^0$ along $\Gamma_1$ and then alternatively computes $u_{1}^{n+1}$ and $u_{2}^{n+1}, ~n=0, 1, 2,...$, as follows:
\beqn
&
\left\{
\begin{aligned}
	\Delta u_1^{n+1} &	= 0, \quad \text{in} ~\Omega_1, \\
	u_1^{n+1} &=   g,  \quad   \text{on} ~ \partial\Omega_1\backslash\Gamma_1 , \\
	u_{1}^{n+1} &=  u_2^{n}, \quad  \text{on} ~\Gamma_1,
\end{aligned}
\right.
\left\{
\begin{aligned}
	\Delta u_2^{n+1} &=  0,  \quad \text{in} ~\Omega_2, \\
	u_2^{n+1} &= g,   \quad \text{on} ~\partial\Omega_2\backslash\Gamma_2 , \\
	u_{2}^{n+1} &=  u_1^{n+1},  \quad \text{on} ~\Gamma_2.
\end{aligned}
\right.
\eeqn

Schwarz presented the earliest alternating Schwarz method to prove Dirichlet's principle.
Later, with the development of parallel computing and computational mathematics, a multitude of interesting DDM have been developed. 
Formally, we can extend the Schwarz method to a general differential operator with boundary condition as follows:
\beq \label{eq:general_equation}	
\mathcal{L} (u) &= f, \quad \text{in}\, \Omega, \\
\mathcal{B} (u) &= g, \quad \text{on}\, \partial\Omega, 
\eeq
where $\mathcal {L} $ can be a negative Laplace operator, $-\Delta $, or any reasonable operator, $\mathcal{B}$ can be $\frac{\partial}{\partial \mathbf{n}} + I $ or any boundary condition, and $f, g$ are two given functions. 
We always presume that the problem \eqref{eq:general_equation} is well-posed. 
Suppose that $\mathcal{D}$ is the operator of an artificial interface transmission condition, 
e.g., Dirichlet, Neumann or Robin. 
We then have a general parallel DDM with domain setting in Figure \ref{fig:diagram_schwarz} listed in Algorithm \ref{alg:general_schwarz}. 
Different from the Schwarz alternating method, the presented Algorithm~\ref{alg:general_schwarz} can be implemented parallel. 

\begin{algorithm}\label{alg:general_schwarz}  
	\caption{A parallel Domain Decomposition Method for Two Subdomains in Figure \ref{fig:diagram_schwarz}}
	\begin{algorithmic}[1]
		\STATE{ Give an initial guess $w_1^{0}, w_2^{0}$ along $\Gamma_1$}
		\FOR{i=1, ...}
		\STATE{Solve for $u_1^{i}$
			\beq \label{alg1:subproblem1}
			\mathcal{L} (u_1^{i}) &= f,\quad \text{in}~ \Omega_1, \\ 
			\mathcal{B}(u_1^{i}) &= g,\quad \text{on}~\partial \Omega_1\backslash \Gamma_1, \\
			\mathcal{D}(u_1^{i}) &= w_1^{i-1},\quad \text{on}~\Gamma_1.
			\eeq}
		\STATE{$w_2^{i} \leftarrow \mathcal{D}\left( u_1^{i}(x_{\Gamma_2}) \right)$}
		\STATE{Solve for $u_2^{i}$
			\beq \label{alg1:subproblem2}
			\mathcal{L}(u_2^{i}) &= f,\quad \text{in}~\Omega_2,  \\
			\mathcal{B}(u_2^{i}) &= g, \quad \text{on}~\partial \Omega_2\backslash \Gamma_2, \\
			\mathcal{D}(u_2^{i}) &= w_2^{i-1},\quad \text{on}~\Gamma_2,
			\eeq}
		\STATE{$w_1^{i} \leftarrow \mathcal{D}\left( u_2^{i}(x_{\Gamma_1}) \right)$}
		\IF{$\Vert w_1^{i} - w_1^{i-1} \Vert / \Vert w_1^{i} \Vert < tol_{\Gamma}$
			\AND $\Vert w_2^{i} - w_2^{i-1} \Vert / \Vert w_2^{i} \Vert < tol_{\Gamma}$
		}
		\STATE{STOP}
		\ENDIF
		\IF{${||u_1^{i} - u_1^{i-1} ||} / {||u_1^{i}||} < tol_{\Omega} $ \AND 
			${||u_2^{i} - u_2^{i-1} ||} / {||u_2^{i}||} < tol_{\Omega} $ }
		\STATE{STOP}
		\ENDIF
		\ENDFOR
	\end{algorithmic}
\end{algorithm}

Here, in Algorithm \ref{alg:general_schwarz}, the stop strategy is the relative error of the current solution with respect to the previous one on the artificial interface or inside the subdomain less than a given tolerance $tol_{\Gamma}$ or $tol_{\Omega}$, respectively. 

\section{Physics-informed neural networks}\label{sec:DeepLearning}
We now introduce the physics-informed neural networks, which are deep fully connected feedforward neural networks, for solving PDE~\cite{raissi2019physics}. 
The entire neural network consists of $L+1$ layers, 
where layer $0$ is the input layer and layer $L$ is the output layer.
Layers $0<l<L$ are the hidden layers.
All of the layers have an activation function, excluding the output layer.
The activation function can take the form of Sigmoid, Tanh (hyperbolic tangents) or ReLU (rectified linear units).

We denote $d_0, d_1, ..., d_L$ as a list of integers, with $d_0, d_L$ representing the lengths of input signal and output signal of the neural network. 
Define function ${\bf T}_l: \mathbb{R}^{d_l} \rightarrow \mathbb{R}^{d_{l+1}}$, $0 \leq l < L$,
\beqn
{\bf T}_{l}({\bf x})={\bf W}_{l} {\bf x} + {\bf b}_{l},
\eeqn
where ${\bf W}_{l} \in \mathbb{R}^{d_{l+1} \times d_{l}}$ and ${\bf b}_{l} \in \mathbb{R}^{d_{l+1}}$. 
Thus, we can simply represent a deep fully connected feedforward neural network using the composite function $h(\cdot;\theta):\mathbb{R}^{d_0} \rightarrow \mathbb{R}^{d_L}$,
\beqn
h(\cdot;\theta) = {\bf T}_{L-1}\circ\sigma\circ {\bf T}_{L-2} \circ \cdots \circ {\bf T}_1 \circ \sigma \circ {\bf T}_0,
\eeqn
where $\sigma$ is the activation function and $\theta:=\{{\bf W}_{l},{\bf b}_{l}\}$ represents the collection of all parameters. 

Solving a general PDE such as \eqref{eq:general_equation} by DNN is a physics-informed minimization problem with the objective $\mathcal{M}(\theta)$ consisting of two terms as follows:
\beq\label{eq:DNN:whole_problem}
\theta^*={\rm argmin}_{\theta}\, \mathcal{M}(\theta):=\mathcal{M}_{\Omega}(\theta)+\mathcal{M}_{\partial \Omega}(\theta), \\
\mathcal{M}_{\Omega}(\theta):=\frac{1}{N_f}\sum_{i=1}^{N_f}
\left\vert \mathcal{L}\left(h({\bf x}_f^{i};\theta)\right)-f({\bf x}_f^i)\right\vert^{2},\\
\mathcal{M}_{\partial \Omega}(\theta):=\frac{1}{N_g}\sum_{i=1}^{N_g} 
\left\vert \mathcal{B}\left(h({\bf x}_g^{i};\theta) \right)-g({\bf x}_g^{i}) \right\vert^{2}, 
\eeq
where $\{{\bf x}_f^i \}_{i=1}^{N_f}$ and $\{{\bf x}_g^i \}_{i=1}^{N_g}$ are the collocation points in the inside and on the boundary, respectively.
The domain term $\mathcal{M}_{\Omega}$ and boundary term $\mathcal{M}_{\partial \Omega}$ enforce the condition that the desired optimized neural network $h(\cdot;\theta^*)$ satisfies $\mathcal{L}(u)=f$ and $\mathcal{B}(u)=g$,  respectively. 

Gradient descent methods can be used to solve this kind of optimization problem; however, from the empirical point of view, more effective and efficient stochastic gradient descents with minibatches are recommended \cite{lecun2012efficient}. 
With regard to the analysis of convergence of stochastic gradient descent, there are many early contributions, such as \cite{bottou-98x, kiwiel2001convergence}.

We should note that the PINN actually provide another kind of discretization scheme for solving PDE. 
As in the cases of FEM or FDM, a basic and interesting question is as follows: is it possible to use $h({\bf x};\theta)$ to approximate the solution of PDE? 
A well-known answer is as follows: if the solution $u({\bf x})$ is bounded and continuous, then 
$h({\bf x};\theta)$ can approximate $u({\bf x})$ to any desired accuracy, given the increasing hidden neurons~\cite{cybenko1989approximation,hornik1991approximation,hornik1989multilayer}. 

\section{Deep domain decomposition methods}
\label{sec:DeepDDM}
Now, it is natural to mix DDM with DNN. 
Here, we provide a general formulation. 
Suppose that we divide the problem \eqref{eq:general_equation} into $S$ subproblems defined in $S$ subdomains as follows:
\beq\label{eq:subproblem}
\mathcal{L}(u_s) &= f_s, \quad \text{in}~\Omega_s, \\
\mathcal{B}(u_s) &= g_s, \quad \text{on}~\partial\Omega_s\backslash\Gamma_s, \\
\mathcal{D}(u_s) &= \mathcal{D}(u_r), \quad \text{on}~\Gamma_s
\eeq
where $f_s$ and $g_s$ are the expressions of $f$ and $g$ in domain $\Omega_s$, respectively, 
and $\Gamma_s$ is the artificial interface of domain $\Omega_s$ to other subdomains. 
To simplify the discussion, we simply denote $u_r$ as the solution of the neighboring subproblems of $u_s$. 
We also denote $\mathcal{D}$ as the operator of the artificial interface transmission conditions.

We denote $h_s({\bf x};\theta)$, $1\leq s \leq S$ as the DNN we used for each subproblem \eqref{eq:subproblem}, i.e., a surrogate of the solution $u_s$,
and denote $\theta_s$ as network parameters of the $k$-th subproblem.
That means different neural networks are used for different subproblems. 
In contrast with the minimization problem \eqref{eq:DNN:whole_problem}, 
we can define subminimization problems as follows:
\beqn
\mathcal{M}_s(\theta; \mathbf{X}_s):=&\mathcal{M}_{\Omega_s}(\theta; {\bf X}_{f_s})+\mathcal{M}_{\partial \Omega_s\backslash\Gamma_s}(\theta;{\bf X}_{g_s})\\
&+ \mathcal{M}_{\Gamma_s}(\theta; {\bf X}_{\Gamma_s}),
\eeqn
where
\beqn
&\mathcal{M}_{\Omega_s}(\theta; {\bf X}_{fs}):=\frac{1}{N_{f_s}}\sum_{i=1}^{N_{f_s}}
\left\vert \mathcal{L}\left(h_s({\bf x}_{f_s}^{i};\theta)\right)-f({\bf x}_{f_s}^i)\right\vert^{2}, \\
&\mathcal{M}_{\partial \Omega_s \backslash\Gamma_s}( \theta; {\bf X}_{g_s}):=\frac{1}{N_{g_s}}\sum_{i=1}^{N_{g_s}} 
\left\vert \mathcal{B}\left(h_s({\bf x}_{g_s}^{i};\theta) \right)-g({\bf x}_{g_s}^{i}) \right\vert^{2}, \\
&\mathcal{M}_{\Gamma_s}(\theta; {\bf X}_{\Gamma_s}):=\frac{1}{N_{\Gamma_s}}\sum_{i=1}^{N_{\Gamma_s}} 
\left\vert \mathcal{D}\left(h_s({\bf x}_{\Gamma_s}^{i};\theta) \right)-
\mathcal{D}\left(h_r({\bf x}_{\Gamma_s}^{i};\theta)\right) \right\vert^{2}, 
\eeqn
with ${\bf X}_{f_s}:=\{{\bf x}_{f_s}^i \}_{i=1}^{N_{f_s}}$,
${\bf X}_{g_s}:= \{{\bf x}_{g_s}^i \}_{i=1}^{N_{g_s}}$ 
and $ {\bf X}_{\Gamma_s}:=\{ {\bf x}_{\Gamma_s}^{i}\}_{i=1}^{N_{\Gamma_s}}$ representing the collocation points in the inside, on the local boundary and on the interface for the $s$-th subproblem, respectively. Denote $\mathbf{X}_s:=\left\{ {\bf X}_{f_s},
{\bf X}_{g_s}, 
{\bf X}_{\Gamma_s}
\right\}$ as the set of whole data.
When the training data $\mathbf{X}_s$ are rearranged into several minibatches at each epoch, we only split the data inside the subdomain.
That means, at each epoch, we first randomly rearrange ${\bf X}_{f_s}$ into $m_s$ disjoint parts $\{ {\bf X}_{f_s}^{1}, {\bf X}_{f_s}^{2},..., {\bf X}_{f_s}^{m_s} \}$, and then let $\mathbf{X}_s^{k}:=\left\{  {\bf X}_{f_s}^{k},
{\bf X}_{g_s}, 
{\bf X}_{\Gamma_s}
\right\}$ be the $k$-th minibatch\footnote{We numerically found that this strategy can make the learning process stable. Here we give a heuristic explanation. Since the number of $N_{g_s}$ or $N_{\Gamma_s}$ is much smaller than that of $N_{f_s}$, once we directly rearrange the total training data, we may have no boundary information in some minibatches. However, for solving a general PDE, the boundary information is quite essential. A more detailed theoretical or numerical analysis would strengthen this point. Considering it's somehow independent with the key inspirit of DeepDDM, we ignore the related numerical discussion.}.
For simplicity, we denote ${\bf W}_s:= \mathcal{D}\left(h_r({\bf X}_{\Gamma_s};\theta) \right) $ as the information that would be transported to the objective subproblem labelled by $s$ from the neighboring subproblems labelled by $r$. Let $\epsilon$ be the learning rate. 
The DeepDDM algorithm for subproblem \eqref{eq:subproblem} is given in Algorithm \ref{alg:DDDM}.

\begin{algorithm*}
	\caption{DeepDDM for the $s$-th subproblem~\eqref{eq:subproblem}, $1\leq s\leq S$.}
	\label{alg:DDDM}
	\begin{algorithmic}[1]
		\STATE{Construct $\mathbf{X}_s$;} 
		\STATE{Initial parameters $\theta^0_s$ and interface information ${\bf W}_{s}^{0}$ along $\Gamma_s$;}
		\FOR{$i=1, 2,...$}\label{step:for_DD}
		\STATE{Set $\theta_s^{i,0}:=\theta_s^{i-1}$;} \COMMENT{Start DDM iteration}
		\FOR{$j=1, 2, ...$}\label{step:for_NN}  
		\STATE{Set $\theta_s^{i,j}:=\theta_s^{i,j-1}$;} \COMMENT{Start DL iteration}
		\STATE{Rearrange randomly training data $\{\mathbf{X}_s^{k} \}_{k=1}^{m_s}$;}
		\FOR{$k=1, 2, ..., m_s$}\label{step:for_minibatch}
		\STATE{
			Update $\theta_s^{i,j}$: 
			$\theta_s^{i,j} \leftarrow \theta_s^{i,j}-\epsilon \nabla_{\theta_s^{i,j}}\mathcal{M}_s(\theta_s^{i,j};\mathbf{X}_s^{k})$; \COMMENT{Update on minibatch}
		}\label{step:argmin}
		\ENDFOR
		\IF{$ \left\vert \mathcal{M}_s(\theta_s^{i,j};\mathbf{X}_s) - \mathcal{M}_{s}(\theta_s^{i,j-\eta};\mathbf{X}_s) \right\vert / \left\vert \mathcal{M}_s(\theta_s^{i,j};\mathbf{X}_s) \right\vert < tol_{\mathcal{M}}$
		}\label{step:eta}
		\STATE{BREAK;}
		\ENDIF
		\ENDFOR \label{step:endfor_NN}
		\STATE{Set $\theta_s^{i}:=\theta_s^{i,j}$;}
		\STATE{Interchange the interface information $\mathbf{W}_s^{i} \leftarrow \mathcal{D}\left( h_r(\mathbf{X}_{\Gamma_s}; \theta_s^{i})\right)$;}
		\IF{$ \left\Vert {\bf W}_s^{i} - {\bf W}_s^{i-1} \right\Vert  \Big{/} {\left\Vert {\bf W}_s^{i} \right\Vert} < tol_{\Gamma}, ~\forall~s $ \label{step:endfor_DDM1}
		}
		\STATE{STOP;}
		\ENDIF
		\IF{$ \left\Vert h_{s} \left( {\bf X}_{f_s}; \theta_s^{i} \right) - h_{s}\left({\bf X}_{f_s};\theta_s^{i-1} \right) \right\Vert / \left\Vert h_{s}\left({\bf X}_{f_s};\theta_s^{i} \right) \right\Vert < tol_{\Omega}, ~ \forall~s $
		}\label{step:endfor_DDM2}
		\STATE{STOP;}
		\ENDIF
		\ENDFOR
	\end{algorithmic}
\end{algorithm*}

We give some remarks as follows:
\begin{itemize}
	\item The loop in Step \ref{step:for_DD} implements the iteration among subproblems, i.e., domain decomposition; the loop in Step \ref{step:for_NN} implements the iteration of the objective function in each subproblem, i.e., deep learning; and the loop in Step \ref{step:for_minibatch} implements the iteration on each minibatch;
	\item The stop criterion for subproblem optimization, shown in Step \ref{step:eta}, is the relative error of current loss and historical loss at step $j-\eta$. Here once a subproblem reach the stop criterion, it will break the loop 5. The interface information will be interchanged after the same loops of all subproblems are broken.
	\item The stop criterion for DDM, shown in Step \ref{step:endfor_DDM1} and Step \ref{step:endfor_DDM2}, is the relative error of solution on interface collocation points and inside collocation points, respectively.
	Here the DeepDDM algorithm will stop after all the subproblems reach stop criteria.
\end{itemize}

\textbf{Computational cost:} 
Obviously, benefiting from DDM, DeepDDM can be implemented in parallel. 
Suppose we fixed the total number of training data, denoted by $N$. 
Once we divide the domain to $S$ subdomains, then the training data for each subproblem is about $N/S$. 
For a fixed network architecture, the training time of one epoch is $O(N/S)$. 
Denote the number of training epoch (DL iteration in Algorithm \ref{alg:DDDM}) for each subproblem as $J_s$. 
Denote the number of DDM iteration as $I_S$ which usually depends on the number of subdomains. 
Then given a fixed network architecture, the total computational cost of DeepDDM is $O({\rm max}_s J_s \cdot I_S \cdot N/S)$. 
For considering the dependency of computation cost on network architecture, it highly depends on the computational device (GPU) 
and the detailed operators in the network. 
For example, if the model size located at a reasonable range and we use NVIDIA V100, for a fixed number of training data, the learning time for one epoch varies little.

\section{Numerical examples}\label{sec:Examples}
In this section, we present a set of systematic numerical results of elliptic problems, including a model problem and an interface problem with a curved interface.
The results of the model problem in different subdomain cases are compared with respect to various network architectures and overlapping size conditions in Section \ref{sec:poisson}.
In Section \ref{sec:discontinuous}, the numerical results of the interface problem with discontinuous coefficients are shown, further demonstrating the effectiveness and robustness of our approach.
Benefiting from the mesh-free strength,
the implementation of domain decomposition can be easily handled even in cases of curved interfaces or more complex geometries.
DeepDDM not only attains promising accuracy but also preserves some properties of DDM existing in conventional discrete strategy.

\textbf{Setting:} In our cases, network architectures are simply set as having equal units for each layer.
The $Layers$ used later indicates the hidden layers in the network, 
and $Units$ indicates the units per layer.
We chose Tanh as activation function and stochastic gradient descent using Adam \cite{kingma2014adam} to update the optimization method.
The mini-batches are set as 64 in this paper, if not otherwise specified.
For the choice of learning rate, the initial learning rates and decay rates change for different situations.
More specifically, learning rates initially range from $10^{-3}$ to $10^{-2}$ and decay every $1-100$ steps with a base of $0.99-0.999$.
The training points herein are all randomly selected. The test points used to calculate errors are regularly sampled by row and column, with 40,000 in all.
The figure of training data indicates the amount of that in the whole domain, and so, the amount of data used in a subdomain is equal to the figure divided by the number of subdomains in the model problem or by the subdomain area ratio in the interface problem. 

\subsection{Model problem}\label{sec:poisson}
The Poisson equation is a simplified form of important equations derived from engineering and physical problems.
The effectiveness of applying algorithms to the Poisson equation is a necessary prerequisite for algorithms applied to more complicated equations.

We consider a Poisson equation with Dirichlet boundary conditions:
\beq\label{eq:poisson}
-\nabla\cdot(\nabla u) & = f,\quad \text{in}~\Omega :=[0, \pi]\times[0, 1],  \\
u(x, y) &= g, \quad \text{on}~\partial\Omega.
\eeq
To obtain an analytic solution, we take, for example, 
$u_{*}=sin(2x)e^{y}$, and substitute it into \eqref{eq:poisson} to compute $f$ and $g$.
We denote the numerical solution by $u_h$ and define the relative $\mathbb{L}^2$ error:
\beqn
\mathcal{E}:=\left( 
\frac{\sum_{i=1}^{n}|u_{*} - u_h |^{2}} 
{\sum_{i=1}^{n}|u_{*}|^{2}}
\right)^{\frac{1}{2}}.
\eeqn

To begin with, we present the results of solving \eqref{eq:poisson} by physics-informed neural networks.
The stop criterion of the inner iteration training process is 
$\left|({\mathcal{M}^{n} - \mathcal{M}^{n-100}})/{\mathcal{M}^{n}}\right|<10^{-3}$, 
or that we exceeded the maximum number of allowed iterations, set as 50,000 here, when the strategy without domain decomposition is used.
Although there is no theory to ensure that the optimization method certainly converges on the global minimum point, our practical experiences indicate that good results would be obtained if the problem is well-posed; the network is large enough; the amount of data is large.

In \autoref{tab:L2error_Net_1D_Dirichlet}, relative $\mathbb{L}^2$ errors between the predicted and the exact solution are presented for different network architectures, while the total numbers of training points are maintained as $N_g=50\times 4$ and $N_f=2500$, respectively.
As expected, if the number of $Layers$ is fixed, more $Units$ corresponds with smaller relative $\mathbb{L}^2$ errors.
Moreover, for a fixed number of $Units$, relative errors would also be reduced when the number of $Layers$ increases.
The general trend suggests that the numerical solutions obtained will be more accurate upon increasing expression capacity of the network.
This observation is in agreement with the results of \cite{raissi2019physics, yang2019adversarial}.

\begin{table*}[!t]
	\centering
	\caption{Relative $\mathbb{L}^2$ error between the predicted and the exact solution for different numbers of hidden layers and different numbers of units per layer.
		Here, the total numbers of training data are fixed to $N_g=50\times 4$ and $N_f=2500$.}
	\label{tab:L2error_Net_1D_Dirichlet}
	\begin{tabular}{c|cccccc}
		\toprule
		\diagbox{$Layers$}{$Units$} 
		& 10 &	20 & 30 & 40  & 50 & 100  \\
		\midrule
		2 & 4.1e-3 & 2.6e-3 & 1.4e-3 & 1.5e-3 & 1.0e-3 & 1.0e-3 \\
		3 & 1.7e-3 & 1.2e-3 & 5.8e-4 & 3.2e-4 & 4.3e-4 & 3.7e-4 \\
		4 & 1.7e-4 & 1.0e-3 & 7.1e-4 & 5.0e-4 & 5.0e-4 & 2.1e-4 \\
		\bottomrule
	\end{tabular}
\end{table*}

\begin{table*}[]
	\centering
	\caption{Relative $\mathbb{L}^2$ error between the predicted and the exact solution for different numbers of training data.
		The network architecture is fixed to $Layers=3$, $Units=50$.}
	\label{tab:L2error_data_1D_Dirichlet}
	\begin{tabular}{c|cccccc}
		\toprule
		\diagbox{$N_g$}{$N_f$} 
		& $100$  & $400$ & $900$ & $1600$  & $2500$ & $10000$  \\
		\midrule
		$10\times 4$ & 7.3e-3 & 2.8e-3 & 3.4e-3 & 6.3e-3 & 9.3e-3 & 2.5e-3 \\
		$30\times 4$ & 4.5e-3 & 2.0e-3 & 1.9e-3 & 1.4e-3 & 1.8e-3 & 1.4e-3 \\
		$100\times 4$ & 3.9e-3 & 9.8e-4 & 7.1e-4 & 6.2e-4 & 4.7e-4 & 1.2e-4 \\
		\bottomrule
	\end{tabular}
\end{table*}

In \autoref{tab:L2error_data_1D_Dirichlet}, we fixed the network architecture, $Layers=3$ and $Units=50$ and recorded the change of relative $\mathbb{L}^2$ error under different numbers of training data.
As expected, as the amount of training data is increased, a more accurate numerical solution is obtained.
Furthermore, the output is more sensitive to $N_g$ than $N_f$.
This observation is in agreement with the results of \cite{raissi2019physics, yang2019adversarial}.

Furthermore, to illustrate the effectiveness of DeepDDM and investigate its properties, 
we present the results of the two-subdomain case and multi-subdomain case.
In practice, we set $tol_{\Omega}=tol_{\Gamma}=10^{-2}$.
The stop criterion of the inner iteration training process is that  
$\left|({\mathcal{M}^{n} - \mathcal{M}^{n-100}})/{\mathcal{M}^{n}}\right|<5*10^{-3}$, 
or that the maximum number of allowed iterations, set as 10,000, is exceeded when DeepDDM is used herein.
Let $\delta$ be overlapping.
For the two-subdomain case, we set $\Omega_1:=[0, \frac{\pi}{2}+\frac{\delta}{2} ]\times[0, 1]$, $\Omega_2:= [\frac{\pi}{2}-\frac{\delta}{2}, \pi ]\times[0, 1] $; however, the numerical solutions are simply defined as follows:
\beqn
u_h = \begin{cases}
	{u}_{h, 1}, ~ x \in(0, \frac{\pi}{2}), \\
	{u}_{h, 2}, ~ x \in(\frac{\pi}{2}, \pi).
\end{cases}
\eeqn
For four-subdomain case, we set
$\Omega_1:=[0, \frac{\pi}{4}+\frac{\delta}{2} ]\times[0, 1]$,
$\Omega_2:= [\frac{\pi}{4}-\frac{\delta}{2}, \frac{\pi}{2}+\frac{\delta}{2} ]\times[0,1]$,
$\Omega_3:= [\frac{\pi}{2}-\frac{\delta}{2}, \frac{3\pi}{4}+\frac{\delta}{2} ]\times[0, 1]$, 
$\Omega_4:=[\frac{3\pi}{4}-\frac{\delta}{2}, \pi]\times[0,1]$,
however the numerical solutions are simply defined as follows,
\beqn
u_h = \begin{cases}
	u_{h,1}, ~ x \in(0, \frac{\pi}{4}), \\
	u_{h, 2}, ~ x \in(\frac{\pi}{4},\frac{\pi}{2}), \\
	u_{h, 3}, ~ x \in(\frac{\pi}{2},\frac{3\pi}{4}), \\
	u_{h, 4}, ~ x \in(\frac{3\pi}{4}, \pi).
\end{cases}
\eeqn
\begin{figure*}[!t]
	\centering
	\includegraphics[width=0.8\textwidth]{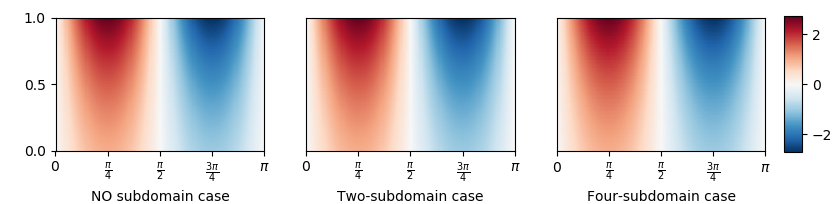}
	\caption{The numerical solutions in the whole domain for different cases. 
		Here, the network architecture is fixed to $Layers=3$, $Units=50$; 
		the total numbers of training data are fixed to $N_g=50\times 4$, $N_{\Gamma}=50$ and $N_f=2500$; and
		$\delta=0.1$ for domain decomposition cases.}
	\label{fig:poisson_numerical_solution}
\end{figure*}

First, for the same situation, same network $Layers=3$, $Units=50$ and same training data $N_g=50\times 4$, $N_{\Gamma}=50$ and $N_f=2500$, the predicted solutions of the Poisson equation for three cases are presented in \autoref{fig:poisson_numerical_solution}.
The results obtained from these cases are close to each other.
\begin{figure*}
	\centering  
	\includegraphics[width=0.45\linewidth]{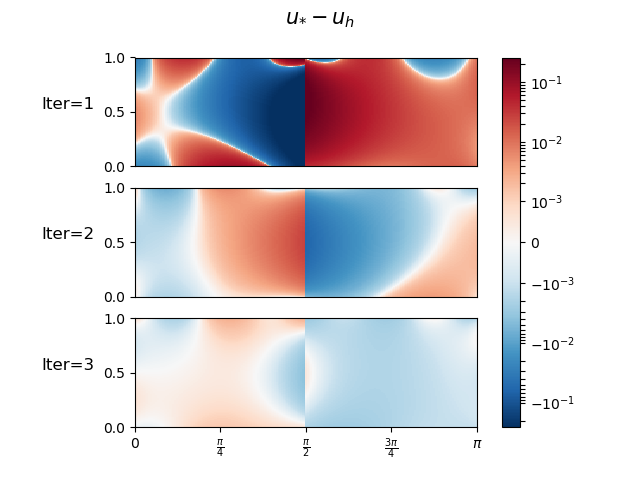}
	\includegraphics[width=0.45\linewidth]{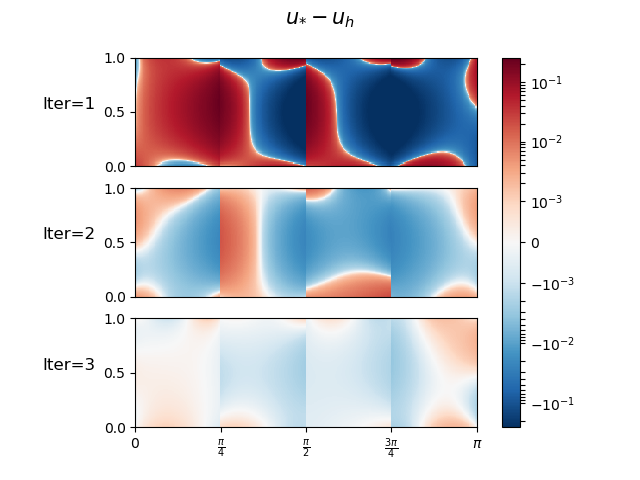}
	\caption{The error $u_* - u_h$ in the whole domain for outer iterations under two subdomains (left) and four subdomains (right). Here, for both cases, the network architecture is fixed to $Layers=3, Units=20$; the total numbers of training data are fixed to $N_g=50\times4$, $N_{\Gamma}=50$ and $N_f=2500$; and $\delta=0.8$.
	}
	\label{fig:poisson_heatmap_DD}
\end{figure*}
In \autoref{fig:poisson_heatmap_DD}, the network is fixed to $Layers=3$, $Units=20$, the numbers of data to $N_g=50\times4$, $N_{\Gamma}=50$ and $N_f=2500$, and overlapping $\delta=0.8$; the errors of solutions in the whole domain are presented for different outer iterations for the two-subdomain case and four-subdomain case.
It is clearly visible that the errors of solutions are relatively larger after the first outer iteration.
However, the numerical solution quickly approximates the exact solution as outer iterations increase.
After the third outer iteration, the algorithm reaches the stop criterion, and the errors of solutions reach the order of $10^{-3}$.
The performance of DeepDDM in the multi-subdomain case is still promising, and the entire convergence process can be clearly observed.

In \autoref{tab:L2error_net_Dirichlet}, the relative $\mathbb{L}^2$ errors and the number of outer iterations are shown for different network architectures for the two-subdomain case and the four-subdomain case.
The accuracy acquired by DeepDDM in \autoref{tab:L2error_net_Dirichlet} is almost the same as that of the results in \autoref{tab:L2error_Net_1D_Dirichlet} in various scales of network architectures.
Otherwise, if we focus on the number of outer iterations, the general trend indicates that the number of outer iterations remain approximately constant no matter how the network architecture changes.
This is consistent with our expectation that the number of outer iterations is independent of the function space, i.e., the network architecture.
Further, all numbers of outer iterations obtained in the four-subdomain case are larger than those of the two-subdomain case.
This observation coincides with the conventional DDM.
As the number of subdomains increases, the number of outer iterations required will increase a little.
The table indicates that DeepDDM not only has the same accuracy as PINN, but also can be implemented in parallel. 
Different training data are employed for different cases, but the impact on our results is extremely limited.%

	\begin{table*}[]
		\centering
		\caption{Relative $\mathbb{L}^2$ error and the number of outer iterations for different network architectures under two-subdomain case and four-subdomain case.
			The figures in parentheses are the numbers of outer iterations.
			For the two-subdomain case, the total numbers of training are fixed to $N_g=100\times 4$, $N_{\Gamma}=100$ and $N_f=10000$.
			For the four-subdomain case, the total numbers of training are fixed to $N_g=50\times 4$, $N_{\Gamma}=50$ and $N_f=2500$.
			Here, $\delta=0.2$ for both cases.}
		\label{tab:L2error_net_Dirichlet}
		\begin{tabular}{c|c|cccccc}
			\toprule
			\makecell{No. of \\subdomains} &
			\diagbox{$Layers$}{$Units$} 
			& 10  &	 20 & 30 & 40  & 50 & 100  \\
			\midrule
			& 2 & 1.9e-3(7) & 4.9e-3(6) & 2.7e-3(7) & 2.6e-3(7) & 4.0e-3(6) & 2.7e-3(7) \\
			Two &
			3 & 3.9e-3(6) & 2.9e-3(7) & 1.0e-3(8) & 2.4e-3(7) & 2.6e-3(7) & 2.5e-3(7) \\
			& 4 & 3.1e-3(7) & 2.6e-3(7) & 2.5e-3(7) & 2.5e-3(7) & 4.3e-3(6) & 4.1e-3(6) \\
			\midrule
			& 2 & 5.4e-3(8) & 6.1e-3(8) & 5.4e-3(8) & 5.6e-3(8) & 5.9e-3(8) & 6.9e-3(8) \\
			Four &
			3 & 3.0e-3(9) & 5.5e-3(8) & 5.8e-3(8) & 5.9e-3(8) & 4.9e-3(8) & 5.7e-3(8) \\
			& 4 & 5.1e-3(8) & 6.1e-3(8) & 5.7e-3(8) & 6.0e-3(8) & 5.7e-3(8) & 5.5e-3(8) \\
			\bottomrule
		\end{tabular}
	\end{table*}
	\begin{figure*}[]
		\centering  
		\includegraphics[width=0.4\linewidth]{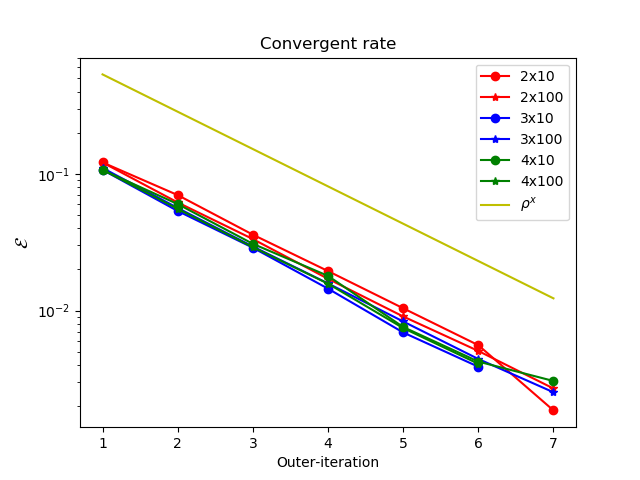}
		\includegraphics[width=0.4\linewidth]{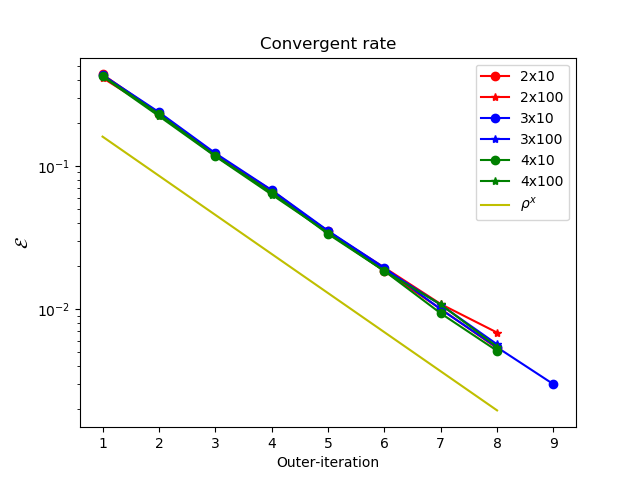}
		\caption{The change in relative $\mathbb{L}^2$ error along with out-iteration in different network architectures under two subdomains (left) and four subdomains (right). Here, the training data are fixed to $N_g=100\times 4$, $N_{\Gamma}=100$ and $N_f=10000$ for the two-subdomain case, the training data are fixed to $N_g=50\times 4$, $N_{\Gamma}=50$ and $N_f=2500$ for the four-subdomain case, and $\delta=0.2$ for both.}
		\label{fig:poisson_convengence}
	\end{figure*}

Further results can be obtained as long as we plot the process of error convergence.
Firstly, if Fourier analysis is used for the domain decomposition model of \eqref{eq:poisson}, we would have the analytic domain decomposition convergence factor:
\beq
\rho:=e^{-k_{min}\delta},
\eeq
where $k_{min}$ is the minimum Fourier frequency, and $\delta$ is the overlap size \cite{gander2006optimized}.
In our case, $k_{min}=\pi$.
The process of error convergence and the analytic convergence factor are shown in \autoref{fig:poisson_convengence}.
Only a few cases are selected; however, other network architecture cases exhibit similar results.
It is clearly observed that the convergence rates of numerical results are always close to that of analysis, regardless of what kind of network architecture is used, in two-subdomain and four-subdomain cases.
This observation coincides with the proposition that the convergence factor depends only on overlap size and domain size.

A more detailed assessment of the effects of overlap sizes is presented in \autoref{tab:out-iter_overlap}.
In particular, we present relative $\mathbb{L}^2$ errors and outer iterations, for both the two-subdomain case and four-subdomain case, with different overlaps from 0.05 to 0.8.
If we fix the overlap and compare results from the two-subdomain case and the four-subdomain case, the latter requires more iterations than the former.
With the increase in the overlapping domain, the outer iterations decline monotonically whether for the two-subdomain case or four-subdomain case.
Moreover, the accuracy would decline as the overlapping domain diminishes, which occurs on account of the decrease of training data in the overlapping domain and the randomness of samples.
It is noteworthy that these observations coincide with the conventional DDM.
\afterpage{
	\begin{table}[!t]
		\centering
		\caption{Relative $\mathbb{L}^2$ error and the number of outer iterations for different sizes of overlap and decomposition.
			Here, the total numbers of training are fixed to $N_g=50$, $N_{\Gamma}=50$ and $N_f=2500$ and the network architecture is fixed to $Layers=3, Units=20$.}
		\label{tab:out-iter_overlap}
		\begin{tabular}{c|cc}
			\toprule
			Overlap & 2x1 decompositions & 4x1 decompositions  \\
			\midrule
			0.05 &  1.6e-2(13)  & 2.1e-2(22)  \\
			0.1  &  9.0e-3(9)   & 1.7e-2(13)  \\
			0.2  &  2.6e-3(7)   & 5.5e-3(8)   \\
			0.4  &  7.6e-4(5)   & 3.3e-3(5)  \\
			0.8  &  9.5e-4(3)   & 6.9e-4(3)  \\
			\bottomrule
		\end{tabular}
	\end{table}
	\begin{figure}[h!]
		\centering
		\includegraphics[width=0.5\textwidth]{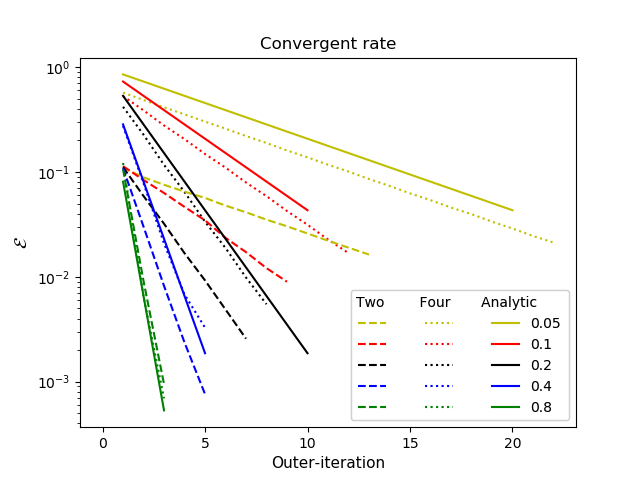}
		\caption{The convergence rates under different overlap sizes and in different subdomain cases.
			Here, the network architecture is fixed to $Layers=3, Units=20$, and the total numbers of training are fixed to $N_g=50\times 4$, $N_{\Gamma}=50$ and $N_f=2500$.
			Dashed line: the two-subdomain case.
			Dotted line: the four-subdomain case.
			Solid line: analytic convergence factors.}
		\label{fig:convergence_overlap}
	\end{figure}
}
Additionally, \autoref{fig:convergence_overlap} presents the comparison of convergence rates for different sizes of overlap under the two-subdomain case and the four-subdomain case.
The line with different colors represents the case with different overlap sizes.
The dashed line and the dotted line respectively represent the two-subdomain case and the four-subdomain case, and the analytic convergence factors are plotted with solid lines.
The figure clearly displays that the numerical convergence rate is closely related to overlap size rather than the number of subdomains.
From the analytical point of view, as long as the overlapping domain expands, the convergence factor would be smaller, and so, the algorithm converges faster in numerical experiments.

The experiments described above demonstrate that DeepDDM effectively inherits the nature of DDM.
Similarly to conventional DDM, as the number of subdomains increases, the necessary iterations also escalate, and the convergence factor shrinks along with the expansion of overlapping domains.
We can perhaps introduce coarse space correction, just like the conventional manner, to improve the performance of DeepDDM.
Moreover, we find that the method of sampling affects the computing results to some extent, and so, the adaptive strategy should be considered and adopted.
However, we leave these questions for future work.

\subsection{Interface problem with discontinuous coefficients}\label{sec:discontinuous}
The example aims to highlight the ability of our approach in handling anisotropic coefficients and a curved interface condition in the governing PDE.
To this end, we consider the case of a steady-state diffusion problem:
\beq \label{eq:discontinuous_diffusion}
-\nabla \cdot (a(\mathbf{x}) \nabla u ) & =  f, \quad \text{in} ~\Omega:=[0, 2]\times[0,2] , \\
\mathcal{B}(u) & =  g, \quad \text{on}~ \partial\Omega, \\
[u]=\left[a(\mathbf{x})\frac{\partial u }{\partial \mathbf{n}}\right] &= 0, \quad \text{on}~ 
\Gamma,
\eeq
where $\Gamma:=\{(x,y)|(x-1)^2+(y-1)^2 =0.25 \}$, $\mathbf{n}$ is the outer normal of the interface and the scalar coefficient $a(\mathbf{x})$ is the piecewise constant function
\beqn
a(\mathbf{x}) = 
\begin{cases}
	1, ~\text{for} ~\mathbf{x} \in \Omega_1, \\
	\alpha, ~\text{for} ~\mathbf{x} \in \Omega_2:=\Omega \backslash\Omega_1, 
\end{cases}
\eeqn
where $\Omega_1:=\{(x,y)| (x-1)^2+(y-1)^2<0.25 \}$. 
The $\alpha$ symbolize the ratio of coefficient.
As $\alpha$ is enlarged, the heterogeneity strengthens, which means that it is more challenging to approximate the exact solution numerically.
In the case of strong heterogeneity in the material, corresponding to a large coefficient ratio, domain decomposition is a recommended approach.
We choose the exact solution as follows,
\beqn
u = \begin{cases}
	\alpha((x-1)^2+(y-1)^2)-0.25(\alpha-1), 
	&~\text{in}~ \Omega_1, \\
	(x-1)^2+(y-1)^2, &~ \text{in}~\Omega_2.
\end{cases}
\eeqn
so that the force term $f$ and boundary condition $g$ in \eqref{eq:discontinuous_diffusion} can be easily obtained.
It is also easy to verify that the above functions exactly satisfy the interface condition in  \eqref{eq:discontinuous_diffusion}.
Here, we choose the boundary operator $\mathcal{B}$ as the identity operator, which means that the Dirichlet boundary condition is employed. 
Using physical coupling conditions between the subdomains, the equation can be written (in a weak sense) in a multidomain formulation,
\beqn
\left\{
\begin{aligned}
	-\Delta u_1 & = f,\quad \text{in}\,\Omega_1, \\
	u_1  &= u_2,\quad \text{on}\,\Gamma,  
\end{aligned}
\right.
\left\{
\begin{aligned}
	-\alpha \Delta u_2  &= f,\quad \text{in}\,\Omega_2, \\
	u &= g,\quad \text{on}\,\partial \Omega_2 \backslash\Gamma, \\
	\alpha \frac{\partial u_2}{\partial \mathbf{n} }  &= \frac{\partial u_1}{\partial \mathbf{n}}, \quad \text{on}\,\Gamma.
\end{aligned}
\right.
\eeqn
The above equations in two subdomains imply that we set $\mathcal{D}$ the identity operator (Dirichlet boundary condition) in $\Omega_1$ and the outer normal operator (Neumann boundary condition) in $\Omega_2$ during applying Alogrithm~\ref{alg:DDDM}. 

To assess the strength of our approach, and by focusing on the curved interface and anisotropic coefficient problems, we present a set of numerical results for various conditions.
The results given demonstrate the effectiveness and the rapid convergence of DeepDDM in various situations.
Moreover, some properties of conventional DDM are also retained in DeepDDM for interface problems, e.g. the iteration is independent of function space; the number of iterations is positively correlated with the ratio of coefficients.
In practice, the tolerance of the stop criterion for DeepDDM is set as $tol_{\Gamma}=tol_{\Omega}=10^{-2}$.
Additionally, the stop criterion of the training process in each subproblem is that
$\left| ({\mathcal{M}^{n} - \mathcal{M}^{n-100}})/{\mathcal{M}^{n}} \right| < 5 \times 10^{-3}$, 
or that the maximum number of allowed iterations, 10,000, is reached.
\begin{figure*}[!t]
	\centering  
	\includegraphics[width=0.5\linewidth]{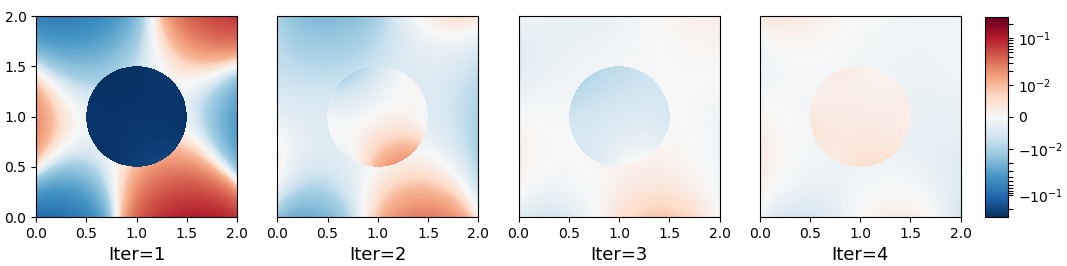}
	\centering
	\includegraphics[width=0.43\linewidth]{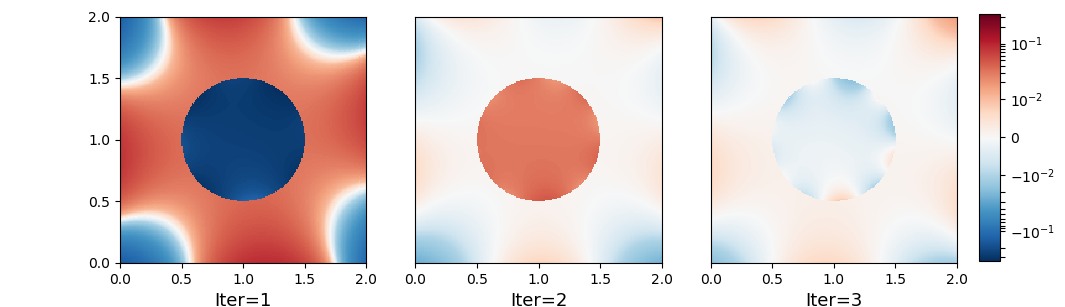}
	\caption{The error $u_* - u_h$ in the whole domain for outer iterations with $\alpha = 2$ (left)  and $\alpha = 20$ (right) in the interface problem. Here, the net architecture is fixed to $Layers=3$ and $Units=20$, and training data are fixed to $N_g=50\times 4$, $N_{\Gamma}=50\times 4$ and $N_f=2500$. 
	}
	\label{fig:discontinuous_heatmap}
\end{figure*}

In the following example, our setup aims to highlight the robustness and effectiveness of the proposed algorithm, even in cases of small network architecture.
The changes in the error of numerical solutions and exact solutions in the whole domain for $\alpha=2$ and $\alpha=20$ are given in \autoref{fig:discontinuous_heatmap}.
Our training dataset consists of $N_g=50\times4$, $N_{\Gamma}=50\times4$ and $N_f=2500$.
All of the training points are randomly sampled.

We chose to represent the latent function $u_h$ using a 3-layer network with 20 units per layer.
As shown by \autoref{fig:discontinuous_heatmap}, the numerical solution quickly approaches the exact solution, along with the subdomain information interchanged on the interface.
Obtained from the test data in the $\mathbb{L}^2$ norm, the relative errors under $\alpha=2$ and $\alpha=20$ are $2.6\cdot 10^{-3}$ and $2.4\cdot 10^{-3} $, respectively.
It is noted that the elementary Dirichlet-Neumann interface condition is used in the strongly discontinuous coefficient problems, implying overlap equal to 0.

The relative $\mathbb{L}^2$ errors and outer-iterations for $\alpha=2$ and $\alpha=20$ are summarized in \autoref{tab:discontinuous_L2error_net}.
Here, the total number of training data is kept fixed at $N_g=50\times4$, $N_{\Gamma}=50\times4$ and $N_f=2500$, respectively.
The key observation here is that as the numbers of $Layers$ and $Units$ increase; i.e., the capacity of the network to approximate more complex functions is enlarged, the numerical accuracy is improved.
When it comes to iterations, they are independent of the network architecture, which is similar to the observation from the previous example.
It is remarkable that the errors of the solution for specific network architecture with $\alpha=20$ are almost equal to the errors of corresponding situations in the $\alpha=2$ case.
Additionally, we also observe that the iteration needed in the $\alpha=20$ case is relatively smaller than that in the $\alpha=2$ case.
In other words, when $\alpha$ becomes larger within a proper range, the number of outer iterations needed is attenuated.
A similar observation has been found in \cite{gander2015optimized}.

	\begin{table*}[]
		\centering
		\caption{Relative $\mathbb{L}^2$ errors and the number of outer iterations for different network architectures.
			These figures in parentheses indicate the number of outer iterations.
			Here, the total number of training data is fixed to
			$N_g=50\times 4$, $N_{\Gamma}=50\times4$ and $N_f=2500$.
		}
		\label{tab:discontinuous_L2error_net}
		\begin{tabular}{c|c|cccccc}
			\toprule
			Cases &
			\diagbox{$Layers$}{$Units$} 
			& 10  &	 20 & 30 & 40  & 50 & 100  \\
			\midrule
			& 2 & 4.4e-3(4) & 3.3e-3(4) & 1.8e-3(5) & 6.1e-4(4) & 2.6e-3(4) & 4.4e-3(3) \\
			$\alpha=2$ & 3 & 6.5e-3(4) & 2.6e-3(4) & 7.7e-3(3) & 1.7e-3(5) & 4.4e-3(4) & 2.0e-3(5) \\
			& 4 & 2.6e-3(5) & 1.9e-3(4) & 2.1e-3(4) & 1.4e-3(5) & 4.0e-3(3) & 3.5e-3(3) \\
			\midrule
			& 2 & 1.2e-2(4) & 1.5e-2(2) & 5.9e-3(3) & 5.5e-3(3) & 3.3e-3(3) & 2.2e-3(3) \\
			$\alpha=20$ & 3 & 1.4e-3(3) & 2.4e-3(3) & 5.1e-3(3) & 3.8e-3(3) & 5.0e-3(2) & 5.3e-3(3) \\
			& 4 & 7.1e-3(3) & 3.9e-3(3) & 5.4e-3(3) & 5.3e-3(2) & 1.4e-3(2) & 2.7e-3(2) \\
			\bottomrule
		\end{tabular}
	\end{table*}
	\begin{figure}[!h]
		\centering
		\includegraphics[width=0.45\textwidth]{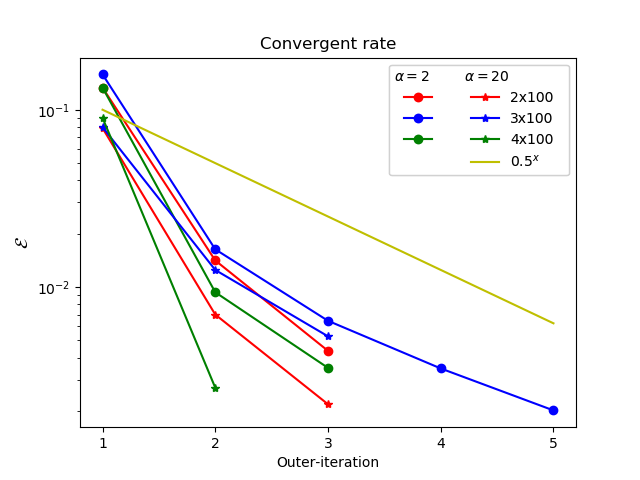}
		\caption{The change in relative $\mathbb{L}^2$ error along with out-iteration in different network architectures and different $\alpha$. Here, the training data are kept fixed to $N_g=50\times 4$, $N_{\Gamma}=50\times4$ and $N_f=2500$.}
		\label{fig:dis_error_iter}
	\end{figure}

Respecting the convergence rate of DeepDDM, \autoref{fig:dis_error_iter} illustrates the change in relative $\mathbb{L}^2$ errors along with outer iterations in various network architectures at $\alpha=2$ and $\alpha=20$. We selected some network architecture cases in \autoref{tab:discontinuous_L2error_net} and then plotted the error as the iterations increase.
According to the yellow line, we can conjecture that the numerical convergence factor is about $0.5$. The theoretical analysis of this interface problem is still an open problem. 


\section{Conclusions}\label{sec:Conclusion}
We have introduced DeepDDM, a novel framework bridging deep learning and domain decomposition, to approximate solutions of PDE according to given equation information.
The presented approach showcases a series of promising results for a model problem and an interface problem.
These numerical results demonstrate that the convergence rate of DeepDDM applied to the Poisson equation is close to the analytical value of DDM.
DeepDDM maintains some properties of conventional DDM: for instance, the outer-iteration depends on the size overlap and the number of subdomains instead of function space.
Furthermore, DeepDDM can easily handle PDE with curve interface and heterogeneity.

This work is seminal in combining deep learning and domain decomposition, and it presents some experiments to provide insights for theoretical study.
The work does not pertain to a stack of fundamental topics behind the approach presented.
Take the following questions as examples.
What is the convergence rate of deep learning solving PDE?
What is the standard to design the best network architecture?
How do we sample the training points for the best result?
Some of these questions are open problems in this field; however, it is still significant to explore the properties of DeepDDM.
Future works should consider applying DeepDDM with the coarse space correction and higher-order interface conditions.

\section*{Acknowledgements}
This research was funded by National Natural Science Foundation of China with grant 11831016 and by Beijing Nova Program of Science and Technology under Grant Z191100001119129. 

\bibliographystyle{ieee_fullname}
\bibliography{yourbibfile}

\end{document}